\documentclass[11pt,a4paper,reqno]{amsart}


\newtheorem{proposition}{Proposition}[section]
\newtheorem{definition}{Definition}[section]
\newtheorem{theorem}{Theorem}[section]

\newtheorem{lemma}{Lemma}[section]

\def\glu{\!\!\!}
\def\hX{\hat{X}}
\def\hT{\hat{T}}
\def\F{\mathcal{F}}
\def\R{\mathcal{R}}
\def\hF{{\hat{\F}}}
\def\hmu{\hat{\mu}}
\def\hf{\hat{f}}
\def\e{\varepsilon}
\def\d{\delta}
\def\htau{\hat{\tau}}
\def\U{{U_r(z)}}
\def\tU{\tau_{\U}}
\def\htU{\hat{\tau}_{\U}}
\def\tG{\tilde{G}}
\def\B{\mathcal{B}}

\def\Z{\mathcal{Z}}
\def\P{\mathcal{P}}
\def\RR{{\mathbb R}}
\def\CC{{\mathbb C}}
\def\NN{{\mathbb N}}
\def\1{{{\mathit 1} \!\!\>\!\! I} }
\renewcommand{\phi}{\varphi}

\renewcommand{\limsup}{\mathop{{\overline {\hbox{{\rm lim}}}}}}

\def\diam{{\hbox{{\rm diam}}}}
\def\path{p}
\def\card{{\hbox{{\rm card}}}}
\def\var{{\hbox{{\rm Var}}}}
\newcommand{\eqdef}{\stackrel{\scriptscriptstyle\rm def}{=}}

\def\Crit{\hbox{{\rm Crit}}}
\def\supp{\hbox{{\rm supp}}}
\def\orb{\hbox{{\rm orb}}}
\def\mod{\hbox{{\rm mod}}}
\def\inte{\hbox{{\rm int}}}

\begin{document}
\bibliographystyle{plain}
\title{Return Time Statistics via Inducing}
\author{H.~Bruin \and B.~Saussol \and S.~Troubetzkoy \and S.~Vaienti}
\thanks{
HB was partly supported by the
PRODYN program of the European Science Foundation.
Also the hospitality of the Centre Physique Th\'eorique,
Luminy, is gratefully acknowledged.\\
BS was supported by FCT's Funding Program and by the Center for Mathematical 
Analysis, Geometry, and Dynamical Systems, Instituto Superior T\'ecnico,
Lisbon, Portugal}
\address{Department of Mathematics,
University of Groningen,
P.O. Box 800, 9700 AV Groningen, the Netherlands}
\email{bruin@math.rug.nl}
%\urladdr{http://www.math.caltech.edu/people/bruin.html}
\address{Laboratoire Ami\'enois de Math\'emathiques Fondamentales et Appliqu\'ees, Universit\'e de Picardie--Jules Verne, 33 rue St Leu,
80039 Amiens Cedex 1, France}
\email{benoit.saussol@mathinfo.u-picardie.fr}
\urladdr{http://www.mathinfo.u-picardie.fr/saussol/}
\address{Centre de Physique Th\'eorique et Institut de Math\'ematiques de
Luminy, CNRS Luminy, Case 907,
F-13288 Marseille Cedex 9, France}
\email{serge@cpt.univ-mrs.fr}
\email{troubetz@iml.univ-mrs.fr}
\urladdr{http://iml.univ-mrs.fr/{\lower.7ex\hbox{\~{}}}troubetz/}
\address{Centre de Physique Th\'eorique, CNRS Luminy, Case 907,
F-13288 Marseille Cedex 9, France}
\email{vaienti@cpt.univ-mrs.fr}
\urladdr{http://www.cpt.univ-mrs.fr/}
\begin{abstract}
We prove that return time statistics of a dynamical system do
not change if one passes to an induced (i.e. first return) map.  
We apply this to show exponential return time statistics
in i) smooth interval maps with nowhere--dense
critical orbits and ii) certain interval maps with neutral fixed points.
The method also applies to iii) certain quadratic maps of the complex plane.
\end{abstract}
\maketitle
\section{Introduction}\label{intro} 
In the last few years the study of {\it return} and {\it hitting} times 
has become an important ingredient for
the statistical characterization of dynamical systems. A 
historical account of this approach can be found in the review
paper~\cite{Coe} or in the introduction of~\cite{hsv}, where 
an extended bibliography is provided. To pose the problem in general terms, 
suppose that $T$ is a measure preserving transformation of a measure space 
$(X,\mu)$ and
$U_x\subset X$ is a neighborhood of a point $x\in X$. 
The two questions which are the fundamental objects of
investigation are:
\begin{enumerate}
\item What is the probability distribution of the first hitting time of 
the set $U_x$ as $\mu (U_x) \to 0$?  

\item What is the probability distribution of the first return time 
for points leaving from $U_x$ as $\mu (U_x) \to 0$?  
\end{enumerate}
For large classes of dynamical systems showing some sort of hyperbolic 
behavior, the answers to those questions are
surprisingly easy and universal. The probability distribution function 
up to a suitable normalization, turns out to be the exponential 1-law
$\exp (-t)$: If $\tau_U(x)$ is the smallest integer $n$ 
such that $T^n(x) \in U$ (so $\tau$ stands for first hitting and for 
first return time), then both
\[
\mu(\{ x \in X; \tau_U(x) \geq \frac{t}{\mu(U)} \}) \to e^{-t} 
\]
and 
\[
\frac1{\mu(U)} \mu(\{ x \in U; \tau_U(x) \geq \frac{t}{\mu(U)} \}) \to e^{-t},
\]
as $\mu(U) \to 0$.
After the works of  Pitskel~\cite{Pit}, Hirata~\cite{hir} and Collet~\cite{col1} 
at the beginning of the nineties which
focused essentially on uniformly hyperbolic dynamical systems, 
new developments have brought at least three improvements:
\begin{itemize}
\item The possibility to treat
larger classes of systems, notably certain non--uniformly hyperbolic systems;

\item The possibility to estimate the error to the asymptotic distribution 
which is closely related to the hyperbolic
character of the transformation and to the influence of sets with short 
recurrence; 

\item The application of
the exponential statistics to the computation of the fluctuations of 
repetition times in the Ornstein-Weiss formula for the
metric entropy.

\end{itemize}
These improvements were the result of four new approaches.
The first approach was originated by Galves and Schmitt~\cite{gs} 
and was formulated in a probabilistic setting for systems satisfying a 
$\phi$-mixing condition. This approach was translated into the dynamical 
systems language by Haydn for
Julia sets~\cite{hay1}, by  Paccaut~\cite{pa} for a large class of non-Markovian maps 
of the interval and by Boubakri~\cite{bo} for some Collet-Eckmann unimodal maps.
Recently Abadi~\cite{ab} extended it for $\alpha$-mixing stationary processes. 
The true probabilistic flavor of this technique led to the discovery of a close 
connection with the Ornstein-Weiss~\cite{ow} theorem on
metric entropy and resulted in a proof of the log-normal fluctuations 
of the repetition times for Gibbs measures in~\cite{cgs} (see~\cite{K} for related results);
see also~\cite{pa} for further developments. 

The second approach is due to Hirata, Saussol and Vaienti~\cite{hsv}.
First for measure-preserving dynamical systems it quantifies the error to the asymptotic 
$\exp(-t)$ distribution giving precise bounds in terms of the mixing properties of the 
systems and of a sharp control of short recurrences (this last point is in itself an interesting
subject suggesting the possibility to formulate a {\it thermodynamics} of return times~\cite{sau,stv}). 
These bounds are then computed for a large class on non-uniformly hyperbolic maps of the interval in~\cite{hsv}. 
We used this method in the present paper to prove the exponential statistics for the class of Rychlik maps
introduced in Section~\ref{Rychlik}. 
A key observation in~\cite{hsv} allows to link the statistics of hitting and return times quoted above, 
namely the distribution of the first return time is close to the exponential law if and only if it is 
close to the distribution of the first hitting time.

The third approach goes back to a probabilistic paper of 
Sevast'yanov~\cite{sev}. This approach was already used by Pitskel~\cite{Pit}
in the context of Markov chains and  Axiom-A diffeomorphisms. The same scheme was applied by Haydn in~\cite{hay1}
and~\cite{hay2} for equilibrium states on Julia sets, 
in the presence of a supremum gap in the 
Perron-Frobenius operator (namely, if $f$ is an H{\" o}lder continuous function on the Julia set of a rational 
map of degree at least~2, then $P(f)> \sup f$, where
$P(f)$ is the topological pressure of $f$). The original technique of 
Sevast'yanov only allowed to prove the convergence to the
exponential 1-law. A recent work by Haydn and Vaienti~\cite{hay3} 
quantifies this technique and thus provides  bounds for the
error estimate in the case of rational maps and parabolic maps of the interval; such bounds improve 
all other existing bounds for the class of maps considered.

The last approach has been introduced by Collet~\cite{col1}; it is an
application of the Gumbel's law for entrance times~\cite{gala}
to a large class of non--uniformly hyperbolic dynamical systems with
exponential decay of correlations for which a tower satisfying Young's
conditions can be constructed~\cite{Y2,Y3}. Collet has proven that the
statistics of closest return to a given point is almost surely
asymptotically Poissonian and gave the fluctuations for the nearest return
to the starting point.

It is useful to remark at this point that the approaches we are describing
also permit to compute the probability distribution of {\it successive}
visits to the set $U_x$: it gives, in all the cases where the exponential
statistics holds, the Poisson distribution $\frac{t^n}{n!}e^{-t}$, where $n$
is the number of visits of $U_x$.

In this paper we propose a new scheme which allows us to compute the
asymptotic distributions for the first and successive return times whenever
these distributions are known on a subset endowed with an induced
structure. 
Its interest lies in the fact
that several non-uniformly hyperbolic systems admit around almost every
point (w.r.t. the invariant measure) a neighborhood where the first return
map acts as the induced transformation and exhibits hyperbolic behavior.
The first application of our approach is for $C^2$ interval maps with
non-positive Schwarzian derivative; we prove the exponential return time
statistics under the additional hypothesis that the map preserves a
conformal measure (see Section~\ref{interval} for details) and that the closure of
the orbit of the critical points has zero measure. Our result 
improves those given by Collet~\cite{col1} and Boubakri~\cite{bo}, in
the sense that, contrarily to them, the growth rate of the derivatives
along the critical orbit plays no role in our theorem. As a second
example, we improve the exponential statistics for the parabolic map
studied in~\cite{hsv}. In both cases, the first return map of the induced
systems belongs to a class of piecewise monotonic maps of the interval
(with countably many pieces), previously investigated by Rychlik; the
induced measure shows exponential decay of correlations with respect to
the Rychlik map and this is the kind of hyperbolic behavior which is
sufficient to prove the exponential statistics for such maps.

In most of the approaches quoted above, the set $U_x$ shrinking to
$\{x\}$ was chosen in the class of cylinders generated by some partition
of the space. The exponential return time statistics around cylinders is what is
needed to compute the fluctuations in the Ornstein-Weiss theorem~\cite{sau}.

Instead of cylinders one can also use balls for the sets $U_x$.
The use of balls has a twofold interest: first it enlarges the class of 
sets where to check the recurrence of rare events, such as entrance time 
(particularly from the perspective of numerical simulations in physical 
situations, see for example~\cite{cs} and~\cite{zt}),
second the distribution of return  times in balls seems to be related 
to other statistical indicators (local dimension and Lyapunov exponents), 
as some work in progress~\cite{stv} suggests.

Balls instead of cylinders were used by Haydn~\cite{hay2} 
and Collet~\cite{col1}.  In this article we prove exponential 
statistics for Rychlik maps around balls. This improves the work 
of~\cite{pa}, who
considered a particular class of Rychlik maps 
(the covering weighted systems introduced in~\cite{lsv1}), and only cylinders
around points. 

The general scheme proposed in this paper could in principle be applied to 
a wide class of dynamical systems for which the induced transformation 
enjoys an exponential
statistics; we think in particular to  rational maps of the Riemann sphere 
(an example is given in Chapter~\ref{complex}) 
and of billiards for which the induced
structures, although not yet studied, seem however more accessible than the original systems. 

As a final observation, we would
like to point out the robustness of the exponential statistics, which persists in 
the non-hyperbolic systems covered in this paper. 
This confirms its central role in the  characterization of recurrence 
for ergodic systems and motivate
further analysis for larger class of transformations especially in dimension greater 
than one and exhibiting singularities.

\section{Statistics via inducing: an abstract theorem}\label{inducing}

Suppose that $(T,X,\mu)$ is an ergodic measure preserving transformation
of a smooth Riemannian manifold $X$.
Let $\hX \subset X$ be an open  set. 
%whose boundary has zero measure.
Furthermore let 
$\hT:\hX \to \hX$ be the first return map. We denote the induced measure by $\hmu$:
the measure preserving transformation $(\hT,\hX,\hmu)$ 
is therefore ergodic.\footnote{We could weaken our 
assumptions by demanding that the map $\hT$ is ergodic instead of 
$T$; the ergodicity of $T$ is however recovered when 
$X = \cup_{n=0}^{\infty} T^{-n}\hX$}
For $z \in X$ we
denote by $\U$ the ball of
radius $r$ centered at $z$ and by $\tU$ (resp.~$\htU$)
the first return time of $\U$ for $T$ (resp.~$\hT$).
For a $\mu$ (resp.~$\hmu$) measurable set $A$ we denote by $\mu_A$ (resp.~$\hmu_A$)
the induced measure on the set $A.$
We suppose that  $(\hT,\hX,\hmu)$ has return time statistics $\hf(t)$:
i.e.~for $\hmu$-a.e.~$z \in \hX, \ \exists \e_z(r) > 0$ with $\e_z(r) \to
0$ as $r \to 0$ such that
\begin{equation}\label{eps_r}
\sup_{t \ge 0} \Big | \hmu_{U_r(z)} \Big ( x \in U_r(z): \htU(x) > 
\frac{t}{\hmu(U_r(z))} \Big ) - \hf(t) \Big | < \e_z(r).
\end{equation} 
Our main theorem of this section is that on $\hX$ the map
$\hT$ obeys the ``same'' statistical
law as $T$:
\begin{theorem}~\label{maintheo}
If $\hf$ is continuous at $0$, then
there exists $f:\RR^+\to [0,1]$ 
such that $N \eqdef \{x: \ f \ne \hf\}$ is countable
and for $\mu$--a.e.~$z \in \hX$ 
and $t\not\in N$, there exists $\d_{z,t}(r) > 0$ with $\d_{z,t}(r) \to
0$ uniformly in $t$ as $r \to 0$ such that
\begin{equation}
\Big | \mu_{U_r(z)} \Big ( x \in U_r(z): \tU(x) > \frac{t}{\mu(U_r(z))} \Big ) - f(t) \Big | < \d_{z,t}(r).
\label{ej}
\end{equation} 
In addition, if $\hf$ is continuous then the convergence is uniform in $t$, which means
there exists $\d_z(r)\to 0$ as $r\to 0$ such that for all $t\in\RR^+$ and $r>0$, $\d_{z,t}(r)\le \d_z(r)$.
\label{thm:1}
\end{theorem}

\underline{Remark:}  the function $\hf$ (resp.~$f$) is  known to be decreasing and
thus continuous everywhere
except for an at most countable set of exceptional points.
The set $N$ is a subset of the points of discontinuity of $\hf.$

\begin{proof}
First of all, if $\mu$ has an atom then $\mu$ is supported on a periodic orbit, 
hence the result is trivial because in this case $f(t)=\hf(t)=1$ for $t\le 1$
and $0$ otherwise. 
We may assume then that $\mu$ has no atoms, 
consequently $\mu(\U)\to 0$ as $r\to 0$ for all $z\in \hX$.

At first we suppose that $z \in \hX$ and assume that $r$ is small
enough that $U_r(z) \subset \hX$. Note that this implies $\hmu_\U=\mu_\U$.

For $x \in \hX$ let $n(x) \eqdef \tau_{\hX}(x)$ be the 
$T$--first return time of $x$ to
$\hX$.  By Kac's Theorem and ergodicity of $\hmu$ we have 
\begin{equation}
A_m(x) \eqdef  \frac{1}{m} \sum_{i=0}^{m-1} n(\hT^ix) 
\mathop{\longrightarrow}_{m \to \infty}
c \eqdef \int_{\hX} n(x) \, d\hmu = \frac{1}{\mu(\hX)}
\label{e3}
\end{equation}
for $\hmu$-a.e.~$x \in \hX.$  Let $G \eqdef \{x \in \hX: \lim_{m \to \infty}
A_m(x) = c\}.$ Clearly $\hmu(G) = 1.$ For all $x \in G$ 
and for all $\e > 0$ there exists $m(x,\e) <
\infty$
such that $|A_m(x) - c| < \e$ for all $m \ge m(x,\e).$
Let $G_m \eqdef \{x \in G: m(x,\e) < m\}$ where we have suppressed the obvious
$\e$ dependence on the set $G_m$.
We choose $M \eqdef M(\e)$ such that $\mu(G_M) > 1 - \e.$

By definition $| \sum_{i=0}^{m-1} n(\hT^ix) - cm | < \e m$
for all $m \ge M$ and all $x \in G_M$.  Thus $\hT^m(x) = T^{cm + s}(x)$
for some $s = s(x)$  
with $|s| < \e m.$
It immediately follows that $\tU(x) = c \htU(x) + s$ with $|s| < \e \htU(x)$
whenever $\htU(x) \ge M$ and $x \in G_M.$

Next we define
\begin{equation}
\tG_M \eqdef \Big \{ z \in G_M: \mu_{U_r(z)}(G_M) 
\mathop{\longrightarrow}_{r \to 0} 1 \Big \}. 
\label{e4}
\end{equation}
By the Lebesgue Density Theorem, $\hmu(\tG_M) = \hmu(G_M)$ and thus $\hmu(\tG_M) > 1 - \e.$
For each $z \in G_M$ there exists $r(z,M,\e) > 0$ such that 
$\mu_{U_r(z)}(G_M) > 1 - \e$ for all $r < r(z,M,\e).$ 
Thus if $R > 0$
is sufficiently small then $\hmu(G_{M,R}) > 1 - \e$,
where $G_{M,r} \eqdef \{z \in \tG_m: r(z,M,\e) > r\}.$

Let $S$ be the set of points $z\in\hX$ for which $\e_z(r)\to 0$ as $r\to 0$.
By assumption $\hmu(S)=1$. Hence $\hmu(G_{M,R}\cap S)=\hmu(G_{M,R})>1-\e$.

Denote 
\begin{eqnarray*}
\F_\U(t) &\eqdef& \mu_{\U}(\{x\in U_r(z): \tau_{\U} > t/\mu(\U)\}), \\
\hF_\U(t) &\eqdef& \hmu_{\U}(\{x\in U_r(z): \htau_{\U} > t/\hmu(\U)\}),
\end{eqnarray*}
and set $B_\U(M)=\{x\in\U : \htau_\U(x) > M\}$.

Assume additionally that $z \in G_{M,R}\cap S$. 
The limiting distribution $\hf(t)$ is continuous at 
$0$ and $\hf(0)=1$, hence if
$r$ is sufficiently small
\[
\mu_{\U}(B_\U(M)^c) = \hmu_{\U}(B_\U(M)^c) \le 1 - \hf(M\hmu(\U)) + \e_z(r) < \e.
\]
Thus for all $t\in [0,\infty)$,
\begin{eqnarray*}
\F_\U(t) &\le& \mu_{\U} \left ( G_{M,R}^c \cup B_\U(M)^c \right ) + \\ 
&\quad& + \mu_{\U}\left ( G_{M,R} \cap B_\U(M) \cap
\left\{x:\tau_{\U}(x) > \frac{t}{\mu(\U)}\right\} \right )  \\ 
&\le& \mu_{\U} (G_{M,R}^c) + \mu_{\U} (B_\U(M)^c) +\\
&\quad& +  \mu_{\U}\left ( G_{M,R} \cap B_\U(M) \cap
\left\{x:\htau_{\U}(x) > \frac{t/(c+\e)}{\mu(\U)}\right\} \right ) \\ 
&\le& 2 \e + \hF_\U\left(\frac{t}{1+\e/c}\right).
\end{eqnarray*}
A similar computation with $1-\F_\U(t)=\mu_\U(\tau_\U\le t/\mu(\U))$ yields 
\begin{equation}\label{eq:***}
\hF_\U\left(\frac{t}{1-\e/c}\right) - 2\e \le \F_\U(t) \leq \hF_\U\left(\frac{t}{1+\e/c}\right) + 2\e.
\end{equation}
For $r$ sufficiently small, $\e_z(r) < \e$. Thus inequality~\eqref{eq:***}
implies:
\begin{equation}
- 3 \e + \hf \left ( \frac{t}{1+\e/c} \right ) 
\le \F_\U(t) \le \hf \left ( \frac{t}{1+\e/c} \right ) + 3 \e .
\label{rr1}
\end{equation}
Define
\begin{equation}
\delta_t(z) \eqdef 
\max \left \{\left |\hf(t) - \hf\left (\frac{t}{1-c^{-1}\e} \right )\right |, 
\left | \hf(t) - \hf \left ( \frac{t}{1+c^{-1}\e} \right ) \right |\right \} .
\label{rr2}
\end{equation}
For each point $t$ of continuity of $\hf$, the function $\delta_t(\e) \to 0$ as $\e \to 0$.
Combining~\eqref{rr1} and~\eqref{rr2} yields:
\begin{equation}
| \F_{U_r(z)}(t) - \hf(t)| \le 3 \e + \delta_t(\e).
\label{eq:****}\end{equation}

We just showed that for any $\e>0$, there exists a 
set $G(\e)$ with $\hmu(G(\e)) > 1-\e$ 
and a real number $R(\e)>0$ such that for all 
$z\in G(\e)$, Inequality~\eqref{eq:****} holds whenever $r<R(\e)$.
The conclusion of the theorem for $\mu$-a.e.~$z \in \hX$  follows from the 
remark that $\mu$-almost every $z\in \hX$ is contained in
the set of full measure $\cap_{n>0} \cup_{m>n} G(1/m)$.

Finally we want to prove the uniform convergence in $t$ in the case $\hf$ is continuous.
Since $\hf(t) \to 0$ as $t \to \infty$, it is 
uniformly continuous; hence $q(\delta) = 
\sup_{0 \le s<t<s+\delta} |\hf(s)-\hf(t)| \to 0$ 
as $\delta \to 0$.

Moreover, $\hf(t)$ is bounded by $1/t$ by Chebychev's inequality.
Hence Equation~\eqref{rr2} gives

\[ 
\delta_t(\e) \le 
\min \left ( q\left (\frac{t\e}{c - \e} \right ) , 1/t \right ) + 3\e.
\]
If $t\e/(c-\e) < \sqrt{\e}$ then $\delta_t(z) \le q(\sqrt{\e}) + 3\e,$ 
while if $t\e/(c-\e) \ge \sqrt{\e}$ 
then $\delta_t(z) \le 1/t \le \sqrt{\e}/(c - \e) + 3\e.$
\end{proof}

\underline{Remark:} Using the same method it is possible to show the counterpart 
of Theorem~\ref{maintheo} for the successive return times and the number of visits.

\section{Piecewise monotonic transformations}\label{Rychlik}

In this section we show that piecewise monotonic maps
of Rychlik's type (even with countably many monotonic pieces) enjoy exponential
statistics. Let $X\subset \RR$ be a compact set and $m$ a Borel regular probability measure on $X$.
The variation of a function $g:X\to\RR$ is defined by
\[
\var\ g \eqdef \sup \{\sum_{j=0}^{k-1} |g(x_{j+1})-g(x_j)| \},
\]
where the supremum is taken along all finite ordered sequences 
$(x_j)_{j=1,\ldots,k}$ with $x_j\in X$. 
The norm $\|g\|_{BV} = \sup|g| + \var\ g$ makes $BV=\{g:X\to\RR:\|g\|_{BV}<\infty\}$
into a Banach space.
We endow $X$ with the induced topology and denote by $\B(X)$ the Borel $\sigma$-algebra of $X$. 
We say that $I\subset X$ is an $X$-interval if there exists some interval $J\subset \RR$ 
such that $I=X \cap J$.

\begin{definition}[$\R$-map]
Let $T:Y\to X$ be a continuous map, $Y\subset X$ open and dense, and  $m(Y)=1$.
Let $S=X\setminus Y$. We call $T$ an {\em $\R$-map} if the following is true:
\begin{enumerate}
\item
There exists a countable family $\Z$ of closed $X$-intervals with disjoint interiors
(in the topology of $X$)
such that $\cup_{Z\in\Z} Z\supset Y$ and for any $Z\in\Z$ the set $Z\cap S$ consists
exactly of the endpoints of $Z$;
\item
For any $Z\in\Z$, $T|_{Z\cap Y}$ admits an extension to a homeomorphism from $Z$ 
to some $X$-interval;
\item
There exists a function $g:X\to [0,\infty)$, with $\var\ g<+\infty$, $g|_S=0$ such that the
operator $P:L^1(m)\to L^1(m)$ defined by
\[
Pf(x)=\sum_{y\in T^{-1}(x)} g(y) f(y)
\]
preserves $m$. In other words, $m(Pf)=m(f)$ for each $f\in L^1(m)$, 
that is $m$ is $g^{-1}$-conformal;
\item
$T$ is expanding:
$\displaystyle \sup_{x\in X} g(x) < 1$.
%For some integer $M\geq 1$ we have 
%$\displaystyle \sup_{x\in X} \left(\prod_{n=1}^M g(T^nx)\right) < 1$.
\end{enumerate}
\end{definition}
We first remark that if $T$ is an $\R$-map, then any iterate $T^n$ is also an $\R$-map
(see Lemma~2 in~\cite{r} and the discussion before for details).
Given a weight function $g$ such that 
$\sup_X g < 1$
%$\sup_X (\prod_{n=0}^M g\circ T^n) <1$ 
and $\var\ g<\infty$, Liverani,
Saussol and Vaienti have shown in~\cite{lsv1} the existence of a conformal measure $m$ which
fulfills the hypotheses, under the additional assumption that $T$ is {\em covering}:
\begin{center}
for any interval $I\subset X$, there exists an integer $N>0$ such
that $$\inf_{x \in Y} P^N\chi_I(x)>0.$$
\end{center}

\begin{proposition}\label{pro:rmap}
Suppose that $X=[0,1]$, $m$ is Lebesgue measure and $T$ is a piecewise monotonic 
transformation on $X$. If
\begin{enumerate}
\item
$T|_{Z\cap X}$ is $C^2$ for each $Z\in\Z$,
\item
$T$ is uniformly hyperbolic, i.e. for some $M>0$, $\inf |(T^M)'| > 1$,
the infimum being taken on the subset of $X$ where $(T^M)'$ is defined,
\item
$T$ has bounded distortion, which means
\begin{equation}\label{distortion}
\sup_{Z\in\Z} \sup_Z \frac{|T''|}{|T'|^2} < \infty
\quad\hbox{ \rm and }\quad
\sum_{Z\in\Z} \sup_Z \frac{1}{|T'|} < \infty,
\end{equation}
\end{enumerate}
then $T$ is an $\R$-map with weight function $g=1/|T'|$ on $Y$ and $g = 0$ on $S$.
\end{proposition}
\begin{proof}
It is enough to check that the variation of $g = 1/|T'|$ is finite.
Since $g$ is $C^1$ on the interior of $Z\in \Z$, the variation is estimated by
\begin{eqnarray*}
\var\ g 
& \leq & 
\sum_{Z\in\Z} \int_Z |g'(t)|dt + 2\sum_{Z\in \Z} \sup_Z g \\
&\leq&
\sum_{Z\in\Z} \int_Z \frac{ |T''(t)| }{ |T'(t)|^2 } dt + 
2\sum_{Z\in \Z} \sup_Z \frac{1}{|T'|}  \\ 
&\leq&
\sup_{Z\in\Z}\sup_Z  \frac{|T''|}{|T'|^2} + 2 \sum_{Z\in\Z} \sup_Z \frac{1}{|T'|} < \infty.
\end{eqnarray*}
\end{proof}

$\R$-maps possess very strong statistical properties, which we collect below
\begin{theorem}[Rychlik~\cite{r}]\label{thm:rychlik}
If $T$ is an $\R$-map then there exists an invariant measure $\mu$ 
which is absolutely continuous with respect to the conformal measure $m$ with 
density $h=\frac{d\mu}{dm}\in BV$. The measures
$m$ and $\mu$ have no atoms. In addition, there exists a partition
(mod $m$) of $X$ into disjoint open sets $X_{i,j}$, 
where $i=1,\ldots,k$ and $j=1,\ldots,L_i$
such that each system $(T^{L_i}|_{X_{i,j}},X_{i,j},\mu_{X_{i,j}})$ is 
mixing and has
exponential decay of correlations for bounded variation observable, 
which means that for some
$C\in(0,\infty)$ and  $\theta\in(0,1)$
\begin{equation}\label{eq:decay}
\left|
\int\glu \phi\circ T^{nL_i} \psi d\mu_{X_{i,j}} - \int\glu\phi d\mu_{X_{i,j}} \int\glu\psi d\mu_{X_{i,j}} 
\right|
\leq
C \|\phi\|_{L^1(m)} (\|\psi\|_{BV})  \theta^n
\end{equation}
for any $\phi\in L^1(m)$ and $\psi\in BV$. 
\end{theorem}
This theorem tells us in particular that $X$ can be 
partitioned (mod $m$) into compact sets
$X_{i,j}$ where some iterate of $T$ is again an $\R$-map with a 
mixing measure equivalent to $m|_{X_{i,j}}$ and which has 
exponential decay of correlations
for bounded variation observable.

\begin{theorem}\label{thm:exp4ryc}
Any $\R$-map $T$ with conformal measure $m$ and invariant mixing measure 
$\mu \ll m$ has exponential return time statistics around balls.
\end{theorem}

\begin{proof}
In the proof we use the estimate given by Lemma~2.4 in~\cite{hsv} to apply Theorem~2.1 in~\cite{hsv}.
We recall the quantities considered there ($N$ is any integer).
\begin{eqnarray*}
a_N(U) &\eqdef& \mu_U(\{x:\tau_U(x)\leq N\}),\\
b_N(U) &\eqdef& \sup \left\{ |\mu_U(T^{-N}V) - \mu(V)| \colon V \hbox{ is measurable} \right\},\\
c(U) &\eqdef& \sup_{k\geq 0} | \mu_U(\{x:\tau_U(x)>k\})-\mu(\{x:\tau_U(x)>k\})|.
\end{eqnarray*}
We denote by $h$ the density of the measure $\mu$ with respect to $m$.
Let $U\subset X$ be an interval in $X$ and set $\tau(U)\eqdef \inf\{\tau_U(x):x\in U\}$.
We now compute an upper bound for $a_N$ and $b_N$.

Obviously $a_N(U) = \sum_{n=\tau(U)}^{N} \mu_U(\{x:\tau_U(x)=n\})$,
and for each $n\geq\tau(U)$
\begin{eqnarray*}
\mu_U(\{x:\tau_U(x)=n\})
&\leq&
\frac{1}{\mu(U)} \int \1_{T^{-n}U}  \1_U d\mu\\
&=&
\int \bigl(\frac{\1_U}{\mu(U)}\bigr)\circ T^n \cdot 
(\1_U-\mu(U)) d\mu + \mu(U)\\
&\leq&
C \frac{m(U)}{\mu(U)} \left\|\1_U-\mu(U)\right\|_{BV} \theta^n + \mu(U).
\end{eqnarray*}
We used for the last inequality the estimate on decay of correlations~\eqref{eq:decay}
given by Theorem~\ref{thm:rychlik}, with $\phi=\1_U/\mu(U)$ and 
$\psi=\1_U-\mu(U)$.
Since $U$ is an interval we have $\|\1_U-\mu(U)\|_{BV} \leq 3$,
hence the summation on $n=\tau(U),\ldots,N$ gives
\begin{equation}\label{eq:an}
a_N(U) \leq \frac{3C}{1-\theta} \frac{m(U)}{\mu(U)} \theta^{\tau(U)} + N\mu(U).
\end{equation}
We consider now $b_N(U)$. The decay of correlations given by~\eqref{eq:decay} 
yields (with $\phi = \1_V$ and $\psi = \1_U /\mu(U)$)
\begin{equation}\label{eq:bn}
b_N(U) \leq \sup_{V\in\B(X)} C \|\1_V\|_{L^1(m)} \|\1_U\|_{BV} \theta^N/\mu(U)
\leq
3C \frac{\theta^N}{\mu(U)}.
\end{equation}
Lemma~2.4 in~\cite{hsv} together with~\eqref{eq:an} and~\eqref{eq:bn}
yield (with $N=\frac{2\log\mu(U)}{\log\theta}$)
\begin{equation}\label{eq:cu}
c(U) \leq c_1 \frac{m(U)}{\mu(U)} \theta^{\tau(U)} + c_2  \mu(U) + c_3 \mu(U)|\log\mu(U)|,
\end{equation}
for some constants $c_1,c_2$ and $c_3$ independent of the interval $U$.

Since $\mu$ has no atoms the countable union $W=\cup_{j=0}^\infty T^{-j}S$ has zero measure.
Hence for $\mu$-a.e. point $z\in X$, the iterates $T^kz$ are well defined and $T^kz\not\in S$.
One easily sees then that $\tau(\U)\to +\infty$ as $r\to 0$ provided $z$ is not periodic. 
Consider the set
\[
G = \{z\in X\setminus W: z \hbox{ is not periodic and } D(z)\eqdef\limsup_{r\to 0} \frac{m(\U)}{\mu(\U)}<+\infty\}.
\]
The Lebesgue Density Theorem tells us that $D(z)$ is $\mu$-almost 
everywhere finite because the 
density $h=\frac{d\mu}{dm}>0$, $\mu$-a.e. 
Since $\mu$ is aperiodic (it is 
ergodic with no atoms), we conclude that $\mu(G)=1$.

Moreover, for all $z\in G$, using~\eqref{eq:cu} we get $c(\U)\to 0 $ as $r\to 0$.
We conclude then by Theorem~2.1 in~\cite{hsv} that for $\mu$-almost all $z\in X$,
\[
\mu_\U\left(\tau_\U > \frac{t}{\mu(\U)} \right) \mathop{\longrightarrow}_{r\to 0} \hbox{{\rm e}}^{-t}
\]
uniformly in $t\in [0,\infty)$.
\end{proof}

\section{Statistic of return times for non-hyperbolic interval maps}\label{interval}

Here we prove that a large class of interval maps enjoy exponential statistics 
for the return time. The strategy is to prove that around almost every point, 
an interval can be found whose first return map is an $\R$-map.

\subsection{Smooth maps with critical points}
Let $T:[0,1] \to [0,1]$ be a $C^2$ interval map with
non-positive Schwarzian derivative.
The Schwarzian derivative
$ST \eqdef \frac{T'''}{T'} - \frac32 ( \frac{T''}{T'})^2$
is a third order differential operator, but it can be replaced by
the condition
\begin{equation}\label{Schw}
\frac{1}{\sqrt{|T'(x)|}} \mbox{ is locally convex.}
\end{equation}
See~\cite{dMvS} for more details. 
Denote the critical set by $\Crit \eqdef \{ x \in [0,1] : T'(x) = 0\}$.

For this class of systems,
the existence of invariant probability measures and their statistical properties
have been frequently studied by means of induced maps (jump transformations
rather than first return maps) and tower constructions.
The situation is best understood for unimodal map with exponential growth
of the derivatives along the orbit of the critical value.
This is known as the {\em Collet-Eckmann condition}, but a slightly stronger
form (involving a slow recurrence condition of the critical point)
has been applied by Benedicks \& Carleson~\cite{BC}.

Every Collet-Eckmann unimodal map has an absolutely continuous invariant 
probability measure $\mu$, $\mu$ has exponential decay of correlations
(w.r.t. to some iterate of the map), satisfies the Central Limit Theorem,
and the Perron-Frobenius operator with suitable weights has a spectral gap~\cite{Y1,KN}. 
An important quantity is the {\em tail behavior} of 
the inducing scheme, i.e. the asymptotic behavior of
$m(\{ x; R(x) > n\})$, where $m$ is the reference measure (Lebesgue) and
$R$ the inducing time (playing the role of the first return time
in our paper). 
Collet-Eckmann maps have exponential tail behavior: 
$m(\{ x; R(x) > n\}) \leq Ce^{-\alpha n}$ for some $\alpha, C > 0$.
In~\cite{BSS} the tail behavior of induced maps over
multimodal maps with arbitrary growth rates on the critical orbits
is computed. The correlation decays were computed
using Young's framework~\cite{Y2}. 

To address the return time statistics of multimodal interval maps,
one can apply Collet's result~\cite{col2}, which yields in our terminology
that if $(T,m)$ has an $\R$-map as induced map, and the tail behavior of
the inducing structure is exponential, then there exists
an invariant probability measure $\mu \ll m$ (Kac's Theorem)
and $(T,\mu)$ has exponential return time statistics.
Boubakri~\cite{bo} studied a unimodal Benedicks-Carleson setting of unimodal
maps in detail, proving Poissonian multiple return time
statistics as well as certain fluctuation results.

In this paper we will assume that the orbit of the critical set is nowhere dense.
This enables us to use a first return map as induced transformation.
On the other hand, we need no information on the tail behavior,
except that $\sum_n m(\{ x ; \tau(x) = n \}) < \infty$, which is equivalent
to the existence of an invariant probability measure $\mu \ll m$.
In particular, the growth rate of the derivatives along the critical
orbit plays no role. Furthermore, we state our results for
multimodal maps and for $|T'|^t$-conformal measures $m_t$.

\begin{theorem}\label{exp4multimodal}
Let $T:[0,1] \to [0,1]$ be as above.
Assume that 
\begin{enumerate}
\item
there exists a non-atomic 
$|T'|^t$-conformal probability measure $m_t$ for some $t > 0$,
\item
$T$ preserves an ergodic probability measure $\mu \ll m_t$, and
\item
$m_t( \overline{ \orb(\Crit) }) = 0$.
\end{enumerate}
Then $(\supp(\mu), T, \mu)$ has exponential return time statistics.
\end{theorem}

\underline{Remark:} A special case is $t=1$, i.e. $m = m_t$ is Lebesgue measure. 
A point $c \in \Crit$ is called {\em non-flat} if there exists
an $\ell$, $1 < \ell < \infty$, such that
$|f(x)-f(c)| = {\mathcal O}(|x-c|^{\ell})$.
If each critical point is non-flat, then
the condition $m( \overline{ \orb(\Crit) }) = 0$ follows from
\[
\orb(\Crit) \mbox{ is } \left\{ \begin{array}{ll}
\mbox{ nowhere dense } &\mbox{ if } \#\Crit = 1,\\
\mbox{ a minimal set } &\mbox{ if } 2 \leq \#\Crit < \infty.
\end{array} \right.
\]  
This was shown by Martens~\cite{martens}
(see also~\cite[Theorem V.1.3']{dMvS}) if $\# \Crit = 1$ 
and by Vargas~\cite{vargas} if $2 \leq \#\Crit < \infty$. 

The case $t < 1$ comes into view when $T$ has a periodic attractor
and the measure $m_t$ is supported on a repelling Cantor set.

\begin{proof}
Let $X = \supp(\mu) \setminus \overline{\orb(\Crit)}$. Obviously $\mu(X) = 1$.
Let $x$ be any recurrent point in $X$; by the Poincar\'e Recurrence Theorem,
this concerns $\mu$-a.e. $x \in X$.

Let $Y$ be the component of $[0,1] \setminus \overline{\orb(\Crit)}$ 
containing $x$ and $n$ be any integer such that
$T^n(x) \in Y$. As $\orb(\Crit) \cap Y = \emptyset$, there exists a neighborhood
$U$ of $x$ such that $T^n$ maps $U$ monotonically onto $Y$.
Let $\hT$ be the first return map to $U$.
By taking $n$ sufficiently large, we can ensure that
$U$ is compactly contained in $Y$ and that $T^i(\partial U) \cap U = \emptyset$
for all $i \geq 0$.

Let $W \subset U$ be any maximal interval on which $\hT|_W$ is monotone.
Since $\orb(\partial U) \cap U = \emptyset$, it follows that
$\hT:W \to U$ is onto, and therefore contains a fixed point, say $p_W$.
The derivative $|\hT'(p_W)|$ is uniformly bounded away from $1$,
see~\cite[page 268]{dMvS}.

Write $\hX = X \cap U$ and $Z = X \cap W$.
Let $\Z$ be the partition of $\hX$ into the sets $Z$.
We will show that $\hT:\cup_{Z \in \Z} Z \to \hX$ is an $\R$-map.

Given intervals $I \subset J$, $J$ is said to contain a 
{\em $\delta$-scaled neighborhood} of $I$ if both components of 
$J \setminus I$ have length $\geq \delta |I|$.
Let $\delta_0 > 0$ be the distance between $\partial Y$ and $\partial X$.
For each branch 
$\hT|_W = T^k|_W$ (where $k \eqdef \tau_{U}(W)$ is the return time), there is an 
interval $W' \supset W$
such that $T^k$ maps $W'$ monotonically onto $Y$.
It follows that $T^k(W')$ contains a $\frac{\delta_0}{|U|}$-scaled neighborhood
of $T^k(W)$. More precisely, if $x,y \in W$, then
$T^k(W')$ contains a $\delta \eqdef \frac{\delta_0}{|T^k(x)-T^k(y)|}$-scaled 
neighborhood of $(T^k(x), T^k(y))$.
This allows us to use the Koebe Principle: 
\begin{proposition}[Koebe Principle]
Let $W' \supset W$ be intervals and assume that
$T^n:W' \to T(W')$ is a $C^2$ diffeomorphism. If $T^n(W)$ contains
a $\delta$-scaled neighborhood of $T^n(W)$, then
there exists $K$ such that the distorsion 
\begin{equation}\label{koebe}
\sup_{s,t \in (x,y) \subset W}
\frac{\hT'(s)}{\hT'(t)} \leq K.
\end{equation}
\end{proposition}

\underline{Remark 1:} If $T$ satisfies~\eqref{Schw}, then 
we can choose $K = K(\delta) = (\frac{1+\delta}{\delta})^2$.
In the general $C^2$ case, we have 
$K = (\frac{1+\delta}{\delta})^2(1+\rho \cdot \sum_{i=0}^n |T^i(W)|)$,
where $\rho$ is a nonnegative number such that $\rho \to 0$ 
as $\sup_{i=0}^n |T^i(W')| \to 0$.
Note that for our purpose, $T^n|_W$ is a first return map. 
Hence $\sum_{i=0}^n |T^i(W)| \leq 1$ for all $W$ and $n = \tau(W)$
the return time of $W$.
For the proof of the proposition and remark,
see~\cite[Chapter III.6]{dMvS}).

\underline{Remark 2:} By taking $U$ small, we can choose $K$ as close to $1$
as we want. It follows that $\inf_{Z \in \Z} \inf_Z |\hT'(x)| 
\geq \frac1K \inf_{Z \in \Z} |\hT'(p_W)| > 1$,
i.e. $\hT$ is expanding.

Since $m_t$ is $g^{-1}$-conformal,
\[
m_t(\hX) = \int_Z |\hT'(x)|^t dm_t \leq \sup_Z |\hT'(x)|^t m_t(Z).
\]
By~\eqref{koebe},
\[
\sup_Z g = \sup_Z \frac1{|\hT'|^t} \leq K_0^t \inf_Z  \frac1{|\hT'|^t} 
\leq K_0^t \frac{m_t(Z)}{m_t(\hX)}.
\]
This shows that $\sum_Z \sup_Z g \leq K_0^t < \infty$.

Let $I \supset J$ be intervals such that $T|I$ is monotone,
and let $L$ and $R$ be the components of $I \setminus J$.
Then $T$ {\em expands the cross ratio} if
\begin{equation}\label{cross ratio}
\frac{|T(I)|}{|T(L)|} \frac{|T(J)|}{|T(R)|}
\geq \frac{|I|}{|L|}\frac{|J|}{|R|}.
\end{equation}
Condition~\eqref{Schw} is equivalent to~\eqref{cross ratio}
(see~\cite[Section IV.1]{dMvS}) which is easily seen to carry over to the 
iterates $T^n$.
It follows that $|\hT'|^{-\frac12}$ is convex on each
branch and therefore $\frac{d}{dx}|\hT'|^{-\frac12} = 
-\frac12 \hT'' |\hT'|^{-\frac32}$ is nondecreasing.
Hence, each branch domain 
$(a_W, b_W) \eqdef W$ contains at most one point $r_W$ at which
$\hT''$ changes sign.
We get
\begin{eqnarray*}
\var_Z g 
&\leq& \var_W g \\
&=& \int_{a_W}^{r_W} \frac{|\hT''(s)|}{|\hT'(s)|^{1+t}} ds 
+  \int_{r_W}^{b_W} \frac{|\hT''(s)|}{|\hT'(s)|^{1+t}} ds  \\
&\leq& \sup_W \frac{1}{|\hT'|^t} \left[ 
\left| \log \frac{\hT'(r_W)}{\hT'(a_W)} \right| + 
\left| \log \frac{\hT'(b_W)}{\hT'(r_W)}\right| \right] \\
&\leq&  2 K_0^t  \log K_0 \, \sup_Z g.
\end{eqnarray*}
Therefore
\begin{eqnarray*}
\var_{\hX} g 
&\leq& \sum_{Z \in \Z} \var_Z g + 2 \sum_{Z \in \Z} \sup_Z g  \\
&\leq& \sum_{Z\in \Z} 2 \left( K_0^t \log K_0 + 1 \right) \sup_Z g
\leq 2 K_0^t \left( K_0^t \log K_0 + 1 \right) < \infty.
\end{eqnarray*}
Thus $\hT$ is a $\R$-map and it satisfies the assumptions of 
Proposition~\ref{pro:rmap}.
By Theorems~\ref{thm:exp4ryc} and~\ref{thm:1}, the conclusion follows 
because the set of $\hX$ satisfying the above assumptions have full measure.
\end{proof}

In the proof of this theorem we even get that $\Z$ is a Markov partition for
the induced map $\hT$ (in fact, each monotone branch is onto). Note that this
is much better than what is necessary, since $\R$-maps need not to satisfy this 
extremely strong topological property.
This supports our belief that this method could be in principle applied to much
more general systems, especially those with singularities.

\subsection{Maps with neutral fixed points}
Let $\alpha\in (0,1)$ and consider the map $T_\alpha$ defined on $X=[0,1]$ by
\[
T_\alpha(x) = \left\{ \begin{array}{ll}
x(1+2^\alpha x^\alpha) & \hbox{{\rm if }} x\in [0,1/2), \\
2x-1 & \hbox{{\rm otherwise.}}
\end{array}
\right. 
\]
Let $\mu_\alpha$ denote the invariant measure absolutely continuous with respect to
Lebesgue (see e.g.~\cite{lsv2} for the existence and properties).
The system $(X,T_\alpha,\mu_\alpha)$ has exponential return 
time statistics around
cylinders of some naturally associated partition~\cite{hsv}; 
Here we prove that this is also true if the neighborhoods are balls.
 
\begin{theorem}
For any $\alpha\in (0,1)$ the system $([0,1],T_\alpha,\mu_\alpha)$ has exponential 
return time statistic.
\end{theorem}

\begin{proof}
Let $\hX_0 \eqdef (\frac12,1]$ and for $n\geq 1$ 
\[
\hX_n \eqdef \{x\in X: T^n(x)>\frac 12\ \hbox{{\rm and for $k = 0,1,\dots,n-1$, }} 
T^k(x)<\frac 12\}.
\]
Fix $n \ge 0$ and let $\hX = \hX_n.$
Let $\hT:\hX\to\hX$ be the first return map to $\hX$.
We then define a partition of $\hX$ by $\Z = \{ Z_p: p=1,2,\ldots\}$, 
where $Z_p=\{x\in \hX: \htau_{\hX}(x) = p\}$. One easily check that $\Z$ is
a partition into intervals, and each branch of $\hT:Z\to\hX$ is monotone and onto. 
Let $m$ be the normalized Lebesgue measure restricted to $\hX$ and let $g(x)=1/|\hT'(x)|$
when $x\in\inte(Z)$ for some $Z\in\Z$ and $g(x)=0$ otherwise.
In the proof of Proposition~3.3 in~\cite{lsv2} it is proved that for some constant $K$
\[
\sup_{x,y\in Z} \frac{\hT'(x)}{\hT'(y)} \leq K
\]
for all $Z\in\Z$. A straightforward computation shows that
\[
\log K \geq \left | \log \frac{\hT'(x)}{\hT'(y)}  \right | \geq 
\left| \int_x^y \frac{\hT''(t)}{\hT'(t)}dt \right|
\]
for all $x,y\in Z$ and $Z\in \Z$. Since $\hT$ is increasing and convex on each $Z\in\Z$,
we get
\[
\int_Z \frac{|\hT''(t)|}{|\hT'(t)|^2}dt \leq \sup_Z g \log K.
\]
Taking into account that $\hT Z=\hX$ for any $Z\in\Z$, we get 
\[
\sup_Z g\leq K\frac{m(Z)}{m(\hT(Z))}\leq K m(Z).
\]
Therefore, since $g$ is $C^1$ in the interior of $Z$ and $g|_{\partial Z}=0$,
\[
\var\ g \leq \sum_{Z\in\Z} \int_Z |g'(t)|dt + 2\sum_{Z\in \Z} \sup_Z g
\leq
K(2+\log K) < \infty.
\]
Finally, it is obvious that $\sup g<1$, hence $(\hX,\hT,m)$ is an $\R$-map.
By Theorems~\ref{thm:exp4ryc}  and~\ref{thm:1} $\mu_\alpha$-almost all points 
inside $\hX$ have exponential return time statistics, and the conclusion follows 
since $\cup_{n\geq 0}\hX_n$ has full Lebesgue measure.
\end{proof}

\section{Complex quadratic maps}\label{complex} 

In this section we apply the main framework to certain polynomials on the 
Riemann sphere $\bar \CC$. 
%This includes hyperbolic and subhyperbolic rational maps
%but also some rational maps with neutral
%periodic points which are widely studied in the literature.
Every rational map $T$ has a $|T'|^t$-conformal measure $m_t$
for some $t \in (0,2]$, see Sullivan~\cite{sul1}. If $T$ is 
hyperbolic on the Julia set $J$, then we can take $t$ equal to 
the Hausdorff dimension 
of $J$, and $m_t$ is equivalent to $t$-dimensional
Hausdorff measure.
In general however, $m_t$ can be supported on a proper subset of $J$;
it can even be atomic. This can be an issue if e.g. $T$ has neutral
periodic points. 
%namely if it is
%supported by a neutral periodic orbit.
%But also in the presence of neutral periodic orbits,
%non-atomic $|T'|^t$-conformal measures can sometimes
%be constructed.

For our results we assume that the orbit of the critical point 
does not densely fill the Julia set, but we do not 
require a supremum gap as in e.g.~\cite{hay1}.

We start with a complex version of $\R$-map and show 
(analogous to Theorem~\ref{thm:exp4ryc}) that under some additional
conditions they have exponential return time statistics.

\begin{definition}[complex Markov maps]
Let $T:Y \to X$ be a continuous map, $Y \subset X$ open subsets of $\bar \CC$, 
and $m_t(Y)=1$ where $m_t$ is a probability measure and $0<t\le 2$. 
Let $S=X\setminus Y$. 
We call $T$ a complex Markov map if the following is true:
\begin{enumerate}
\item
There exists a countable family $\Z$ of pairwise disjoint open discs
such that $\cup_{Z\in\Z} Z = Y$;
\item
For any $Z\in \Z$, $T:Z \to X$ is a conformal diffeomorphism (with $T(Z)=X$), and with bounded distortion:
\begin{equation}\label{compdist}
\sup_{Z\in\Z} \sup_Z \frac{|T''|}{|T'|^2} < \infty;
\end{equation}
\item
The measure $m_t$ is $|T'|^t$-conformal;
\item
The map $T$ is expanding: $\displaystyle \inf_{y \in Y} |T'(y)| > 1$.
\item 
Let  $\Z_k = \Z \vee T^{-1} \Z \vee \cdots \vee T^{-(k-1)} \Z$
is the $k$-th join of the partition $\Z$.
The domains $Z \in \Z_k$ are uniformly {\em convex-like}, by which we mean
\begin{equation}\label{path}
\sup_k \sup_{Z \in \Z_k} \sup_{x \neq y \in Z} 
\frac{\path_Z(x,y)}{|x-y|} < \infty,
\end{equation}
where $\path_Z(x,y)$ denotes infimum of the lengths of the path 
in $Z$ connecting $x$ and $y$.
\end{enumerate}
\end{definition}

Analogous to Theorem~\ref{thm:rychlik} one can show that complex Markov maps 
have an invariant measure with exponential decay of correlations.
\begin{theorem}\label{thm:wal}
Let $T$ be a complex Markov map as above.
There exists an invariant ergodic probability measure $\mu$ equivalent to $m_t$.
Moreover $(X,T,\mu)$ is mixing with exponential decay of correlations~:
There exists $C>0$ and $\rho\in(0,1)$ such that for any $f$ Lipschitz and $g\in L_\mu^1$
\begin{equation}\label{doccomp}
\left|\int f \cdot g\circ T^n d\mu -\int fd\mu \int g d\mu\right| 
\le C \|f\|_{Lips} \|g\|_{L_\mu^1}\rho^n,
\end{equation}
where $\|f\|_{Lips} = \sup |f| + \sup_{x\neq y} \frac{|f(x)-f(y)|}{|x-y|}$.
\end{theorem}

\begin{proof}
First let us make the following remark.
As $T$ is expanding,
a straightforward calculation shows that \eqref{compdist}
implies that
\[
K := \sup_{n \in \NN} \sup_{Z \in \Z_n} \sup_{w \in Z}
\frac{|(T^n)''|}{|(T^n)'|^2} < \infty
\]
and $\sum_{Z_n \in\Z_n } \sup_Z |(T^n)'|^{-t} < \infty$
uniformly in $n$.
Let $P_t$ be the Perron-Frobenius operator
\[
(P_t f)(z) = \sum_{ T(y) = z } \frac{f(y)}{|T'(y)|^t}.
\]
Let us write $\| f \|_s$ for the seminorm
$\sup_{z \neq z'} \frac{|f(z)-f(z')|}{|z-z'|}$.
If $\| f \|_{Lips} < \infty$, we have
\begin{eqnarray*}
\| P_t^nf\|_s &\leq&
\sup_{z,z' \in Y} \frac{1}{|z-z'|}
|P_t^n(f)(z) - P_t^n(f)(z')| \\
 &\leq&
\sup_{z,z' \in Y} \frac{1}{|z-z'|}
\left\{
\sum_y
\frac{ |f(y)-f(y')| }{ |(T^n)'(y)|^t } + \right. \\
&&\hskip 1cm
\left. \sum_{y'} |f(y')|
\left| \frac{1}{|(T^n)'(y)|^t} -  \frac{1}{|(T^n)'(y')|^t} \right|
\right\}
\end{eqnarray*}
Here we summed over the pairs $y,y'$ in the same atom $Z_n \in \Z_n$
with $T^n(y) = z$ and $T^n(y') = z'$.
As $T:Z \to X$ is onto for every $Z \in \Z$, these pairs are well-defined.
The first term in the above expression is bounded by
\[
\sum_{Z_n} \sup_{y,y' \in Z_n}
\frac{|f(y)- f(y')|}{|y-y'|} \frac{1}{|(T^n)'(y)|^{1+t}}
\leq \| f \|_s \| P_t^n \1\|_{\infty}
\sup_{y \in Y} \frac{1}{|(T^n)'(y)|}.
\]
Next we use the Mean Value Theorem and \eqref{path} to estimate
\begin{eqnarray*}
\left| \frac{1}{|(T^n)'(y)|^t} -  \frac{1}{|(T^n)'(y')|^t} \right|
&\leq&
\tilde K \cdot \left( \frac{1}{|(T^n)'(w)|^t} \right)' |y-y'| \\
&=&
\frac{\tilde K |t| \, |y-y'|}{|(T^n)'(w)|^{t-1}}. \frac{|(T^n)''(w)|}{|(T^n)'(w)|^2},
\end{eqnarray*}
for some $w$. The constant $\tilde K$ is an upper bound in
\eqref{path}. This gives for the second term
\begin{eqnarray*}
&& \sup_{Z \in \Z} \sup_{z,z' \in Z}
\sum_{y'} \frac{ f(y') }{ |z-z'| }
\left| \frac{1}{|(T^n)'(y)|^t} -  \frac{1}{|(T^n)'(y')|^t} \right| \\
&& \hskip 1cm \leq
\| f \|_{\infty} \sum_{Z_n} \sup_{y,y' \in Z_n}
\frac{1}{|y-y'| \cdot |(T^n)'(y)|}
\left| \frac{1}{|(T^n)'(y)|^t} -  \frac{1}{|(T^n)'(y')|^t} \right| \\
&& \hskip 1cm \leq
\| f \|_{\infty}
\sum_{Z_n} \sup_{w \in Z_n} \frac{\tilde K |t|}{|(T^n)'(w)|^t}
\frac{|(T^n)''(w)|}{|(T^n)'(w)|^2} \\
&& \hskip 1cm \leq
2 K \tilde K \| f \|_{\infty} \| P_t^n\1 \|_{\infty}.
\end{eqnarray*}
Let $\theta = \sup |T'|^{-1} \in (0,1)$.
Then
\[
\frac{ \| P_t^nf \|_s }{ \| P_t^n\1 \|_{\infty} }
\leq
\theta^n \| f \|_s + 2 K \tilde K \| f \|_{\infty}.
\]
Obviously 
$\| P^n_tf\|_{\infty} \leq \| f \|_{\infty} \| P^n\1\|_{\infty}$.
Therefore also\[
\frac{ \| P_t^nf \|_{Lips} }{ \| P_t^n\1 \|_{\infty} }
\leq
\theta^n \| f \|_{Lips} + (2 K \tilde K + 1) \| f \|_{\infty}.
\]
By construction
$\| P_t^n \1 \|_{\infty} = \sum_{Z\in\Z} \sup_Z \frac{1}{|(T^n)'|^t} < \infty$
uniformly in $n$.
This allows us to use the Tulcea-Ionescu \& Marinescu theorem,
which shows that $P_t$ is a quasicompact operator.
Since each branch is onto, $1$ is a simple eigenvalue with Lipschitz 
eigenvector $h>0$, and is the unique eigenvalue on the unit circle.
Consequently, for any Lipschitz function $f$ we have
\[
\|P_t^n(fh)-h \int\glu f d\mu\|_{Lips}\le C_0\|fh\|_{Lips}\rho^n,
\]
for some $C_0>0$ and $\rho\in(0,1)$.
It follows that the correlations between Lipschitz functions
$f$ and $L_\mu^1$ functions decays exponentially fast:
if $g \in L^1_\mu$ we find
\begin{eqnarray*}
\left|\int f \cdot g\circ T^n d\mu - \int\glu fd\mu \int g d\mu\right| 
&=&
\left|\int (P^n_t(fh)-h\int\glu f d\mu) \cdot g dm_t \right|\\
&\le& \|P^n_t(fh)-h \int \glu f d\mu\|_\infty \, \int |g| dm_t \\
&\le& C \|f\|_{Lips} \|g\|_{L_\mu^1}\rho^n,
\end{eqnarray*}
for $C = C_0 \|h\|_{Lips} \|1/h\|_\infty$.
\end{proof}
\begin{theorem}\label{thm:exp4compryc}
Any complex Markov map $T$ as defined above with $|T'|^t$-conformal 
measure $m_t$ for some $t>1$
admits an invariant mixing measure $\mu$ equivalent to $m_t$ with exponential 
return time statistics around balls.
\end{theorem}

The proof of Theorem~\ref{thm:exp4compryc} is similar to that of Theorem~\ref{thm:exp4ryc}.
The additional difficulty is that the decay of correlation is given for Lipschitz function
and not characteristic function of balls. The usual way to overcome this problem is to 
approximate balls by an union of small cylinders. This is the object of the
next lemmas.
Note that if we restrict ourselves to cylinder sets, then the theorem
is valid for all $t > 0$.

As before, let $\Z_k = \Z \vee T^{-1}\Z \vee\cdots\vee T^{-k+1}\Z$ 
denote the $k$-dynamical partition
and $\B_k$ the $\sigma$-algebra generated by $\Z_k$.

\begin{lemma}\label{lem:phimixing}
Let $T$ be a complex Markov map.
There exists some constants $\Gamma$ and $\gamma>0$ such that for any 
set $A\in\B_k$ and Borel set $B$ we have 
\begin{equation}\label{eq:mixing}
\left| \mu( A \cap T^{-n-k} B ) - \mu(A) \mu(B) \right| 
\le \Gamma \mu(A)\mu(B) \exp(-\gamma n).
\end{equation}
\end{lemma}
\begin{proof}
We can rewrite~\eqref{eq:mixing} as
\[
\left| \int h^{-1} \cdot P_t^k(h\1_A) \cdot \1_B\circ T^n d\mu - 
\mu(A)\mu(B) \right| \le \Gamma \mu(A)\mu(B) \exp(-\gamma n).
\]
where $h=\frac{d\mu}{dm_t}$ is Lipschitz and bounded away from 
below (see~\cite{wal}), and
$P_t$ is the Perron-Frobenius operator as defined in the previous proof.
Recall also that
$K \eqdef \sup_{k\in\NN} \sup_{Z\in \Z_k} \sup_Z 
\frac{|(T^k)''|}{|(T^k)'|^2} <\infty$.
Since $T^k:A\to X$ is one-to-one and onto when $A\in\Z_k$, the above facts imply that 
\[
\sup_{k\in\NN}\sup_{A\in\Z_k} \frac{\| h^{-1} P_t^k(h\1_A) \|_{Lips}}{\mu(A)} <\infty.
\]
The lemma follows now from~\eqref{doccomp} in Theorem~\ref{thm:wal}
by taking $f = P_t^k(\1_A)/\mu(A)$ and $g = \mu(A) \1_A$.
\end{proof}

Given $S\subset \CC$ we denote by $\B_k(S)$ the smallest element of the $\sigma$-algebra
$\B_k$ containing $S\cap J$.
\begin{lemma}\label{lem:ball}
There exists $\alpha>0$ such that for any $z\in Y$, $r>0$ and $k>0$ we have
\begin{equation}\label{eq:ball}
\mu\left(\B_k(U_r(z))\cap \B_k(U_r(z)^c) \right) \le \exp(-\alpha k),
\end{equation}
where $U_r(z)$ denotes the ball of radius $r$ about the point $z$.
\end{lemma}
\begin{proof}
Let us assume for simplicity that $\diam(Y)\le 1$.
We denote by $\inte_J S$ the interior of a subset $S\subset J$ in the 
induced topology.
Since $T$ is expanding, there exists $a \in (0,1)$ such that $\diam (Z)< a^k$ 
for any $Z\in \Z_k$ and integer $k$.
By the Markov property, the bounded distortion and the conformality 
of the map $T$ we can find some constant $c>0$ such that the following property
holds:
For any integer $k$ and cylinder $Z\in\Z_k$ there exists $p_Z\in \inte_J Z$ and $r_Z>c\cdot\diam(Z)$
such that $U_{r_Z}(p_Z) \cap \inte_J Z' = \emptyset$ for any 
$Z \neq Z'\in\Z_k$ different from $Z$.
Given $Z\neq Z' \in\Z_k$, we have $d(p_Z,p_{Z'})>\max(r_Z,r_{Z'})$; 
in particular, $U_{r_Z/2}(p_Z)\cap U_{r'_Z/2}(p_{Z'})=\emptyset$.

Let $x\in Y$, $r>0$ and consider the following partition
\[
\P=\{Z\in\Z_k: Z\subset \B_k(U_r(z))\}.
\]
Let $\P_n =\{Z\in\P: a^n <\diam (Z) \le a^{n-1}\}$.
Any $Z\in\P_n$ is a subset of the 
annulus $S_n(z,r) \eqdef U_{r+a^n}(z) \setminus U_{r-a^n}(z)$.
Thus there exists $\card (\P_n)$ disjoint balls or radius at least $a^n/2$
inside $S_n(z,r)$.  Since the area of $S_n(z,r)$ is equal to
$4\pi r a^n$ when $a^n\le r$ we get
\[
\card (\P_n) \le \frac{4\pi r a^n}{c^2 a^{2n}/4} =\frac{8\pi}{c^2} r a^{-n} \le \frac{8\pi}{c^2}a^{-n}.
\]
Obviously when $a^n>r$ we also have $\card(\P_n) \le \frac{8\pi}{c^2} \le \frac{8\pi}{c^2} a^{-n}$.

Since the measure $m_t$ is $|T'|^t$-conformal and the map itself is conformal 
and Markov we have for some constant $c_1$, $m_t(Z) < c_1 \diam(Z)^t$ for any $Z\in\Z_k$ for some $k$.
The previous inequalities imply (observe that $\P_n=\emptyset$ if $n\le k$
and recall that $t > 1$) 
\[
\begin{split}
\mu(\bigcup_{Z\in\P}Z) &= \sum_{n>k} \mu\Big(\bigcup_{Z\in \P_n}Z\Big)\\
&\le \sum_{n>k} \max\{\mu(Z):Z\in\P_n\} \cdot \card(\P_n)\\
&\le \sum_{n>k} c_1a^{-t} a^{tn} \frac{8\pi}{c^2} a^{-n} \\
&= \frac{8\pi c_1}{(a-a^t)c^2} a^{(t-1)k}.
\end{split}
\]
Taking $\alpha\in(0,(1-t)\log a)$ sufficiently small gives the result.
\end{proof}

\begin{proof}[Proof of Theorem~\ref{thm:exp4compryc}]
The proof closely follows the one of Theorem~\ref{thm:exp4ryc}, we use 
the same notation $a_N$ and $b_N$.
Let $U=U_r(z)$ and $k,N\in\NN$ to be chosen later on. 
By Lemma~\ref{lem:phimixing} and Lemma~\ref{lem:ball}
\[
\begin{split}
a_N(U) 
&\le \sum_{n=\tau(U)}^{N} \frac{1}{\mu(U)}\mu(\B_n(U)\cap T^{-n}U)\\
&\le (1+\Gamma) \sum_{n=\tau(U)}^N \mu(\B_n(U))\\
&\le \frac{1+\Gamma}{1-\exp(-\alpha)} \left[ N\mu(U) + \exp(-\alpha \tau(U))\right].
\end{split}
\]
Similarly, we get by Lemmas~\ref{lem:phimixing} and 
(twice) Lemma~\ref{lem:ball}
\[
\begin{split}
b_N(U) 
&\le \Gamma \mu(\B_k(U)) \exp[-\gamma(N-k)] +
 \frac{1}{\mu(U)} \mu(\B_k(U)\cap \B_k(U^c)) \\
&\le \Gamma \exp[-\gamma (N-k)] + 
\frac{1}{\mu(U)} \left( \Gamma \exp[-\alpha k - \gamma(N-k)] +
\exp[-\alpha k] \right),
\end{split}
\]
for all $k \leq N$.
Taking $k = -\frac{2}{\alpha}\log\mu(U)$ and $N=2k$ gives
\[
b_N(U) \leq \Gamma \mu(U)^{\frac{2\gamma}{\alpha}} + 
\Gamma \mu(U)^{1+\frac{2\gamma}{\alpha}} + \mu(U).
\]
For all non periodic points  $z\not\in \cup_{k\in\NN} \partial\Z_k$ we have $\tau(U_r(z))\to\infty$,
which implies $c(U_r(z))\to 0$ as $r\to 0$. Since this concerns $\mu$-almost 
all points in $J$, the theorem is proved.
\end{proof}

We will apply these results to quadratic maps on $\CC$.
Induced systems have been used before for rational maps, notably
by Aaronson et al. \cite{adu}. They consider parabolic rational
maps (i.e. rational maps whose Julia sets contain no critical point
but rationally indifferent periodic points),
and establish the existence of an invariant measure $\mu \ll m_t$,
where $m_t$ is a $t$-conformal measure with $t = \mbox{HD}(J)$.
Moreover, $\mu$ is finite if and only if
$t \cdot \min_p \frac{a(p)+1}{a(p)} > 2$.
Here the minimum is taken over all parabolic points $p$ and $a(p)$ is
such that $T^q(z) = z + \alpha(z-p)^{a(p)} + \dots$ for the appropriate
iterate $q$. It follows that $\mu$ is finite only if $t > 1$,
which is the hypothesis in Theorem~\ref{thm:exp4compryc}.
It is to be expected therefore that parabolic rational maps with a finite
invariant measure $\mu \ll m_t$ have exponential return time statistics
on balls.

\begin{theorem}\label{thm:quadratic}
Let $T(z) = z^2 + c$ be a quadratic map on $\CC$ such that $T$ is not infinitely
renormalizable (see discussion below) and its Julia set $J$ contains no
Cremer points.
Suppose that also for some $t>1$
\begin{itemize}
\item $J$ supports a
non-atomic $|T'|^t$-conformal measure $m_t$,
\item $m_t(\overline{\orb(\Crit)}) = 0$, and
\item $T$ preserves a probability measure $\mu$ equivalent to $m_t$. 
\end{itemize}
Then $(\supp(\mu), T, \mu)$ has exponential return time statistics 
on disks.
\end{theorem}

The Hausdorff dimension of the Julia set $\mbox{HD}(J) > 1$ for
each parameter $c$ in the Mandelbrot set $\mathcal M$
(see Zdunik, \cite{zd});
the two exceptions $c = 2$ and $c=0$ are easy real one-dimensional
cases.
So assuming that these parameter values allow a $t$-conformal
measure with $t = \mbox{HD}(J)$, the condition
$t > 1$ is no restriction. By continuity, there are many parameters  
close to $\mathcal M$ that give $t > 1$ as well.

Note that we allow $T$ to have parabolic points or Siegel disks.
The map $T$ is called {\em renormalizable} if there exist
open disks $W_0$ and $W_1$, $0 \in W_0 \subset W_1$, such that $T^n:W_0 \to W_1$
is a two-fold covering map for some $n \geq 2$ and $T^{in}(0) \in W_1$ for
all $i \geq 1$. If there are infinitely many integers $n$ such that this is possible,
$T$ is {\em infinitely renormalizable}.
We assume that $T$ is not infinitely renormalizable and has no Cremer point
to be able to use Yoccoz' puzzle construction. In particular the result
that for each $z \in J$, the puzzle pieces $P_n(z)$ containing 
$z$ shrink to $z$ as $n \to \infty$ is important. For an exposition of Yoccoz' 
puzzles and the proof of these statements we refer to~\cite{Mil}.

Yoccoz' results have been extended to certain infinitely renormalizable
polynomials by Lyubich and Levin \& van Strien,~\cite{Lyu,LvS}.
For reasons of simplicity, we have not tried to extend Theorem~\ref{thm:quadratic}
to these cases; we prefer to work with a single set of initial puzzle 
pieces $P_0$.
Prado~\cite{prado} has shown that in all of the above cases, the conformal measure $m_t$
is ergodic.

\begin{proof}
Let $z \in \supp(\mu) \setminus \overline{\orb(\Crit)}$
be arbitrary, and let $U \supset V \owns z$ be open disks 
such that $U \cap \overline{\orb(\Crit)} = \emptyset$ and $V$ is 
compactly contained in $U$. Let $\log \delta = \mod(U \setminus V) > 0$
be the modulus of $U \setminus V$.
Assume also that $V$ is {\em convex-like} in the sense that
$\sup_{x \neq y \in V} p_V(x,y)/|x-y| < \infty$, where $p_V$ is
as in \eqref{path}.

The strategy is to find a subset $\hX$ of $V$ such that 
$T^n(\partial \hX) \cap \hX = \emptyset$ for all $n \geq 0$, and 
then we can invoke Theorem~\ref{thm:exp4compryc}.

If $J$ is a Cantor set, we can assume that $\partial V$ is contained in
the Fatou set $F$, which is the basin of $\infty$ in this case.
Moreover, there are no neutral or stable periodic orbits.
Thus each point in $F$, in particular $\partial V$, converges to $\infty$. 
It follows that $T^n(\partial V)$
intersects $V$ for at most finitely many $n \geq 0$.
Let $\hX$ be the component of 
$V \setminus \cup_{n \geq 0} T^n(\partial V)$
containing $z$. Then $T^n(\partial \hX) \cap \hX = \emptyset$ for 
all $n \geq 0$ as required.
Since $\hX$ consists of the intersection and difference of
at most finitely many, convex-like disks, we find
$\sup_{x \neq y \in \hX} p_{\hX}(x,y)/|x-y| \leq C(\hX) < \infty$.

Assume now that $J$ is connected and by Yoccoz' results also locally connected.
Let $F_i$, $i \geq 0$, be the periodic components of the Fatou set, with
$F_0$ the basin of $\infty$. 
Since $T$ is a polynomial (with exceptional point $\infty$),
$J = \partial F_0$.
There are only finitely many such components,
and by Sullivan's Theorem
(see e.g.~\cite{sul2}), every $z' \in F$ is eventually mapped into $\cup_i F_i$.
We construct a special forward invariant subset $G$ of $\overline{\cup_i F_i}$.

\begin{enumerate}
\item
Consider the renormalization $T^n:W_0 \to W_1$ of the highest possible period $n$.
(If $T$ is not renormalizable, then we just take $T:\CC \to \CC$.)
It is known that $W_1$ contains an $n$-periodic point $p$ with at least two external
rays, say $A_0$ and $A'_0$.
The existence of such external rays (and the arcs $A_i$ defined later on)
is guaranteed by results initiated by Douady, 
see~\cite{Pet} and references therein.)
Let $G_0 = \{ z \in \CC; |z| > 10 \}$.

\item
If $F_i$, $i \geq 1$, contains a stable periodic point, let $G_i$ be
a disk compactly contained in $F_i$ such that $T^{per(F_i)}(G_i) \subset G_i$
and $\orb(\Crit \cap F_i) \subset G_i$.
There is at least one $per(F_i)$-periodic point $p_i$ in the boundary
of $F_i$. Let $A_i$ be a smooth compact arc connecting 
$p_i$ and $G_i$ such that
$T^{per(F_i)}(A_i) \subset A_i \cup G_i$.
Since $p_i \in \partial F_0$, there is also an external ray $A'_i$ landing at
$p_i$.

\item
If $F_i$ contains a parabolic point $p_i$ in its boundary
such that each $z \in F_i$ converges to $p_i$, let $G_i \subset F_i$
be a disk such that $T^{per(F_i)}(G_i) \subset G_i$,
$\orb(\Crit) \cap F_i \subset G_i$, and 
$\partial G_i \cap \partial F_i = \{ p_i \}$.
Let $A'_i \subset F_0$ be an external ray landing at $p_i$.
\end{enumerate}
Let 
\[
G  = \bigcup_{i \geq 0} (G_i \cup \{ p_i \} \cup A_i \cup A'_i).
\]
Then $G$ is connected, and for some $N \in \NN$, $\cup_{j = 1}^N T^j G$ is forward invariant.
We start a Yoccoz puzzle construction by putting, for $n \geq 0$,
\[ 
P_n = \{ \hbox{components of } \bar \CC \setminus T^{-n} G \}
\] 
For each $n \geq 1$, $T$ maps any element of of $P_n$ into an element of
$P_{n-1}$. Using the arguments in~\cite{Mil}, one can show that the diameters of 
the elements $Y \in P_i$ tend
to $0$ as $i \to \infty$, unless $Y$ eventually intersects a Siegel
disk.

We can assume that the point $z \notin \cup_n T^{-n}(G)$.
Moreover, since $\orb(\Crit)$ densely fills
the boundary of any Siegel disk, and $m_t(\overline{\orb(\Crit)}) = 0$,
we can assume that $z$ does not lie on the boundary of a Siegel disk.
Find $n$ so large that the element $Y$ of $P_n$ containing
$z$ is contained in $V$.
Let $\hX = V \cap Y$. Then
$T^n(\partial \hX) \cap \hX = \emptyset$ for all $n \geq 0$.
Note also that $T^n(Y)$ is bounded by finitely many smooth curves
of $\partial G_i$, $A_i$ and $A'_i$. At worst these curves end in a
logarithmic spiral, namely as they approach the hyperbolic periodic
points $p_i$. Therefore also $T^n(\hX)$ is convex-like, and obviously
simply connected. It follows that $\sup_{x \neq y \in \hX} p_{\hX}(x,y)/|x-y|
\leq C(\hX) < \infty$.

The rest of the argument works for both $J$ locally connected and 
$J$ a Cantor set.
Let $\hT:\hX \to \hX$ be the first return map to $\hX$.
Then $\hT$ is defined on a countable collection $\Z$ of
disks $Z$. 
The modulus $\mod(U \setminus \hX) \geq \log \delta$,
and for each branch $\hT = T^{\tau}:Z \to \hX$ there exists
a disk $Z' \supset Z$ such that $T^{\tau}$ maps $Z'$ univalently onto
$U$.
It follows from the Koebe $\frac14$-Theorem (e.g.~\cite[Theorem 1.4]{CG})
that the distortion of $\hT|_Z$ is uniformly bounded:
\begin{equation}\label{dist2}
\sup_{x,y \in Z} \frac{|\hT'(x)|}{|\hT'(y)|} \leq K = K(\delta)
\end{equation}
More precisely, take $x \in Z$ and let $[x,y] \subset Z$ be a straight
arc containing $x$ such that
$\hT''/\hT'$ varies little on $A$.
Then
\[
\frac{ |\hT''(x)| }{ |\hT'(x)| } 
\leq 
\frac{2}{|y-x|} \left| \int_x^y \frac{ \hT''(u) }{ \hT(u) } du \right|
\leq
\frac{2}{|y-x|} \log \frac{|\hT'(y)|}{|\hT(x)|}.
\]
To estimate this, let $B$ be the maximal round disk centered at $x$ contained
in $Z'$. Let $\delta_0$ be the radius of $B$; we have 
$\delta_0 |\hT'(x)| = {\mathcal O}(\delta)$.
Define $f(w) = \frac{\hT(x+\delta_0 w) - \hT(x)}{\delta_0 \hT'(x)}$.
Then $f$ is a univalent map on the unit disk with $f'(0) = 1$. By
Theorem~1.6. of~\cite{CG}, we obtain 
\begin{eqnarray*}
\frac{|\hT'(y)|}{|\hT'(x)|^2} 
&=& \frac{1}{|\hT'(x)|} \log \frac{|f'(\frac{|y-x|}{\delta_0})|}{|f'(0)|} \\
&\leq& \frac{1}{|\hT'(x)|} \left[
\log(1+\frac{|y-x|}{\delta_0}) - 3\log(1-\frac{|y-x|}{\delta_0}) \right] \\
&\leq& \frac{5|y-x|}{\delta_0 |\hT'(x)|},
\end{eqnarray*}
provided $y$ is close to $x$.
Therefore
$\frac{ |\hT''(x)| }{ |\hT'(x)|^2 } \leq \frac{10}{\delta_0 |\hT'(x)|}$,
proving that 
\begin{equation}\label{dist3}
\sup_{Z \in \Z} \sup_Z \frac{|\hT''|}{|\hT'|^2} < \infty.
\end{equation}
The distortion bound $K(\delta)$ applies by the same argument also
to iterates $\hT^n$ of the induced map.
Therefore each domain $Z \in \Z_n$ is not much less convex than $\hX$:
\[
\sup_k \sup_{Z \in \Z_k} \sup_{x \neq y \in Z} \frac{p_Z(x,y)}{|x-y|} \leq C
\]
for some $C = C(\delta,\hX) < \infty$.

We can also assume that $V$ and hence $\hX$ is so small that 
$\inf_{\cup Z} |\hT'(x)| > 1$.
Thus $\hT$ is hyperbolic.
Recall that $m_t$ is the $|T'|^t$-conformal measure of $T$.
As $\hT$ is a first return map, it is straightforward to show
that $\hat m_t = \frac{m_t}{m_t(\hX)}$ is a $|\hT'|^t$-conformal
map for $\hT$. In fact, $\hat m_t$ can also be constructed
using Sullivan's techniques~\cite{sul1}.
%By~\eqref{dist2}
%\begin{equation}\label{dist4}
%\sum_{Z\in \Z} \sup_Z \frac{1}{|\hT'|^t} \leq K \sum_{Z\in\Z} \frac{m_t(Z)}{m_t(\hX)}
%= \frac{K}{m_t(\hX)} < \infty.
%\end{equation}
Now we can invoke Theorem~\ref{thm:exp4compryc} with the invariant 
measure $\hat \mu = \frac{1}{\mu(\hX)} \mu$.
\end{proof}

\end{document}